\documentclass[11pt,fleqn,a4paper]{amsart}
\usepackage{a4wide}
\usepackage[]{graphicx}
\usepackage{psfrag}

\def\cal{\mathcal}
\newcommand{\wt}[1]{\widetilde{#1}}
\newcommand{\sier}[1]{{\cal O}_{#1}}
\newcommand{\C}{{\mathbb C}}
\renewcommand{\P}{{\mathbb P}}
\newcommand{\Q}{{\mathbb Q}}
\newcommand{\Z}{{\mathbb Z}}
\newcommand{\ep}{\varepsilon}
\newcommand{\vp}{\varphi}
\newcommand{\la}{\lambda}
\newcommand{\al}[1]{\alpha_{#1}}
\newcommand{\lra}{\longrightarrow}

\newcommand{\Rank}{\mathop{\rm Rank}}
\newcommand{\Pic}{\mathop{\rm  Pic}}

\newtheorem{propo}{Proposition}
\newtheorem{conj}{Conjecture}
\newtheorem*{question}{Question}
\newtheorem*{propstar}{Proposition}

\begin{document}

\keywords{Surface singularity, superisolated singularity, universal
       abelian cover, geometric genus.}
\subjclass[msc2000]{32S25}

\title{Universal abelian covers of superisolated singularities}

\author[Jan Stevens]{Jan Stevens}
\email {\sf stevens@math.chalmers.se}
\address
     {Matematiska Vetenskaper,
G\"oteborgs universitet,  SE 412 96 G\"oteborg, Sweden.} 

\begin{abstract}
We give explicit examples of Gorenstein surface singularities with
integral homology sphere link, which are not complete intersections.
Their existence was shown by  Luengo-Velasco, Melle-Hern\'andez  and
N\'emethi, thereby providing counterexamples to the universal abelian
  covering conjecture of Neumann and Wahl.
\end{abstract}
\maketitle

The topology of a normal surface singularity does not determine the
analytical invariants of its equisingularity class. Recent partial
results indicated that this nevertheless could be true under two
restrictions, a topological one, that the link of the singularity is
a rational homology sphere, and an analytical one, that the
singularity is $\Q$-Gorenstein. Neumann and Wahl conjectured that
the singularity is then an abelian quotient of a complete
intersection singularity, whose equations are determined in a simple
way from the resolution graph \cite{nw}. Counterexamples were found
by Luengo-Velasco, Melle-Hern\'andez  and N\'emethi \cite{lmn}, but
they did not compute the universal abelian cover of the
singularities in question. The purpose of this paper is to provide
explicit examples.

We give examples both of universal abelian covers and of Gorenstein
singularities with integral homology sphere link, which are not
complete intersections. We mention here one example with both
properties, see Propositions \ref{cuspprop} and \ref{verderprop}.
\begin{propstar} The Gorenstein singularity in $(\C^6,0)$ with
ideal generated by the maximal minors of the matrix
\[
\begin{pmatrix} u & y-z & 8s^2-9y-9z & w \\ 5s^2-9y-9z & v & w &
u^2-46s^3+54ys+54zs+us^2
\end{pmatrix}
\]
and three further polynomials
\begin{gather*}
  8u^2s+(8s^2-9y-9z)(4s^2-9y-9z)+27(y-z)^2\;,\\
  8(5s^2-9y-9z)us+w(4s^2-9y-9z)+27(y-z)v\;, \\
  8(5s^2-9y-9z)^2s+(u^2-46s^3+54ys+54zs+us^2)(4s^2-9y-9z)+27v^2\;,
\end{gather*}
has integral homology sphere link. The exceptional divisor on the
minimal resolution is a three-cuspidal rational curve with
self-intersection $-1$.\\
The singularity is the universal abelian cover of a hypersurface
singularity with equation  $$
27t^2+(4s^3+u^2)^3+2u^5s^2-20u^3s^5-4us^8+u^4s^4=0\;, $$ which is
the suspension of an irreducible plane curve singularity with
Puiseux pairs $(3,2)$ and $(10,3)$. The link has as first homology
group $\Z_3$ and the  resolution graph is
\[
\unitlength=30pt
\def\ci{\circle*{0.26}}
\def\vi{\ci}
\def\mbt#1{\makebox(0,0)[t]{$\scriptstyle #1$}}
\begin{picture}(3.5,2.5)(-1.25,-1.25)
\put(-1,-1){\line(1,1){.9}}  \put(1.5,0){\line(-1,0){1.37}}
\put(-1,-1){\vi} \put(-1,-1.2){\mbt{-3}} \put(1.5,0){\ci}
\put(1.5,-.2){\mbt{-7}}
\put(-1,1){\line(1,-1){.9}} \put(-1,1){\ci}
\put(2.5,1.2){\makebox(0,0)[b]{$\scriptstyle -3$}}
\put(0,0){\circle{0.26}} \put(0,-.2){\mbt{-1}}
\put(2.5,1){\line(-1,-1){1}}   \put(2.5,1){\vi}
\put(2.5,-1.2){\mbt{-3}} \put(2.5,-1){\line(-1,1){1}}
\put(2.5,-1){\vi}
\end{picture}
\]
\end{propstar} \noindent This graph occurs in \cite[Examples 2]{nw}
as an example violating the so-called semigroup condition. It is not
the graph of a superisolated singularity, nor is the universal
abelian cover of the type which can occur as cover of such a
singularity. But our method of computing for superisolated
singularities also yields the stated equations.

A surface singularity $f\colon (\C^3,0) \to (\C,0)$ is called
superisolated if blowing up the origin once resolves the
singularity. Writing the  equation $f=f_d+f_{d+1}+\cdots$ in its
homogeneous parts the condition becomes that $f_d$ defines a plane
curve with isolated singularities, through which the curve
$f_{d+1}=0$ does not pass. These singularities were introduced by
Luengo to give counterexamples to the smoothness of the $\mu$-const
stratum \cite{lu}. They provide a class of singularities lying
outside the usual examples (they are degenerate for their Newton
diagram), yet they are easy to handle because most computations
reduce to the plane curve $f_d=0$ alone. To obtain a singularity
with rational homology sphere link from an irreducible curve one
needs that it is also locally irreducible, so it is a rational curve
with only (higher) cusp singularities. The universal abelian cover
of such a superisolated singularity is a cyclic cover of degree $d$.

The typical example of a superisolated singularity has the form
$f_d+l^{d+1}$, where $l$ is a general linear function. This suggests
to generalise the problem and study coverings of the whole Yomdin
series $f_d+l^{d+k}$. In general these singularities  are not
quasi-homogeneous, which complicates the computations. But the two
limiting cases of the series are homogeneous: for $k=\infty$ we have
the non-isolated homogeneous singularity $f_d$, while $k=0$ gives
the cone over a smooth plane curve of degree $d$. To find a cyclic
cover of degree $d$ we note that in both cases we have an
irreducible curve $C$, embedded by a very ample line bundle $H$ of
degree $d$.
 The affine cone over $C$ has local ring $\bigoplus
H^0(C,nH)$. Let $L$ be a line bundle with $dL=H$. The ring
$\bigoplus H^0(C,nL)$ is the local ring of a quasi-homogeneous
singularity (the quasi-cone $X(C,L)$), which is a $d$-fold cover of
the cone over $C$. For a smooth curve there are $d^{2g}$ line
bundles $L$ satisfying $dL=H$, but for singular curves the number is
less. In particular, in the case of interest to us, where $C$ is a
rational cuspidal curve, the Jacobian of $C$ is a unipotent group
and there is a unique line bundle $L$ satisfying $dL=H$. For this
line bundle we determine the ring $\bigoplus H^0(C,nL)$. We find in
this way the $d$-fold cover of the singularity $f_d$. This is again
a non-isolated singularity. If we understand its equations well
enough, we can write down a series of singularities, and quotients
of suitable elements are coverings of singularities in the series of
$f_d$. It is also worthwhile to look at other elements of the
series; in this way the equations in the proposition above were
obtained.
We succeed for the case $d=4$. But for a number of curves of higher
arithmetical genus, where we have  determined the structure of the
ring $\bigoplus H^0(C,nL)$, the equations of the non-isolated
singularity are complicated and have no apparent structure.
Therefore we do not write them, nor do we give  a series of
singularities.

Superisolated singularities with $d=5$ give counterexamples
\cite{lmn} to the conjectures by N\'emethi and Nicolaescu \cite{nn}
on the geometric genus of Gorenstein singularities with rational
homology sphere link. We have made a detailed study of this case. If
the splice quotient of Neumann and Wahl exists, and that is the case
for an exceptional curve with at most two cusps, it has the
predicted geometric genus.

In the first section we discuss the relation between topological and
analytical invariants. The next section recalls the conjectures of
N\'emethi and Nicolaescu, and of Neumann and Wahl. The third section
contains generalities on superisolated singularities, while the next
one presents the results for the case $d=4$. The section on degree
five focuses on the geometric genus. The next section treats
generalisations. In the last section we discuss the conjectures in
the light of our examples.

The computer computations for this paper were done with the computer
algebra system {\sc Singular} \cite{GPS}. Thanks are due to the
referee for his detailed comments and suggestions.

\section{Analytical  and topological invariants}

The conjectures on $\Q$-Gorenstein surface singularities with
rational homology sphere link fit in the wider framework of the
following question.

\begin{question}
Which discrete data are needed to know a normal surface singularity?
\end{question}

One interpretation of `knowing a singularity' is that we can write
down equations. As we only have discrete data, such equations
necessarily describe a family of  equisingular surfaces. At the very
least one should know the geometric genus of the singularity, the
most basic analytical invariant, which is constant in equisingular
families (this requirement excludes the multiplicity). The first
thing needed is the topology of the singularity, or what amounts to
the same, the resolution graph.

Let $\pi\colon (\wt X,E)\to (X,p)$ be any resolution of a normal
surface singularity $(X,p)$ with exceptional set $E=\bigcup_{i=1}^r
E_i$. The intersection form $E_i\cdot E_j$ is a negative definite
quadratic form. It can be coded by a weighted graph, with vertices
$v_i$ corresponding to the irreducible components $E_i$ of the
exceptional set with weight $-b_i=E_i\cdot E_i$; two vertices are
joined by an edge, if the corresponding components intersect, with
edge weight  the intersection number. The numerical class $K$ of the
canonical bundle of $M$ provides a characteristic vector: by the
adjunction formula we have $E_i\cdot(E_i+K)=2p_a(E_i)-2$. The
resolution graph of $\pi$ is the weighted graph obtained  by giving
the vertices as additional weight the arithmetical genus of the
corresponding component, written in square brackets. It encodes all
necessary information for the calculus of cycles on $E$, but it
carries no information on the singularities of $E$. Therefore one
usually requires that all irreducible components are smooth and that
they intersect transversally; this is called a good resolution. All
edge weights are then equal to 1, and are not written. According to
custom we also suppress the vertex weights $[0]$ and $-2$. By the
resolution graph of a singularity we mean the graph of the minimal
good resolution.

For rational singularities  the resolution graph suffices to `know'
the singularity, by a celebrated result of Artin. 
By Laufer 
and Reid 
the same holds for minimally elliptic singularities. Beyond these
classes of singularities this is no longer true, as shown by the
following example of Laufer, see \cite[Example 6.3]{pi}; it occurs
on several places in \cite{nem}, starting with 2.23.

\subsection{Example} Consider the resolution graph:
\[
\unitlength=30pt
\def\ci{\circle*{0.26}}
\def\mbt#1{\makebox(0,0)[t]{$\scriptstyle #1$}}
\begin{picture}(2,1)(-.26,-.5)
\put(2,0){\line(-1,0){1.87}} \put(1,0){\ci}
\put(0,.3){\makebox(0,0)[b]{$\scriptstyle -1$}} \put(2,0){\ci}
\put(0,0){\circle{0.26}} \put(0,-.3){\mbt{[1]}}
\end{picture}
\]
The singularity is numerically Gorenstein, with $K=-3E_0-2E_1-E_2$
on the minimal resolution, where $E_0$ is the elliptic curve, which
is intersected by the rational curve $E_1$ in the point $P$, and
$E_2$ is the remaining rational curve. Let $\sier{E_0}(-Q)$ be the
normal bundle of $E_0$ on the minimal resolution. One has
$E_0\cdot(E_0+K)=2Q-2P$, as divisor on $E_0$. So by adjunction
$2P=2Q$, if the singularity is Gorenstein.

To compute $p_g=H^1(\sier{})$ one can use computation sequences for
S.S.-T. Yau's elliptic sequence; a convenient reference is
\cite{nem}. Starting from $Z=E_0+E_1+E_2$ one considers in order the
cycles $Z+E_0$, $Z_1=Z+E_0+E_1$ and $Z_2=-K=Z_1+E_0$. One looks at
the long exact cohomology sequences of the short exact sequences
connecting these cycles. The most relevant ones are
\[
0 \lra {\cal O} (-Z-E_0) \lra{\cal O}(-Z) \lra  \sier {E_0}(Q-P)
\lra 0
\]
and
\[
0 \lra {\cal O} (-Z_1-E_0) \lra {\cal O}(-Z_1) \lra  \sier
{E_0}(2Q-2P) \lra 0\;.
\]
As $H^1(\sier{}(-Z_1-E_0))=H^1(\sier{}(K))=0$ one finds $p_g\leq 3$.
Moreover, if $2P\neq 2Q$, the singularity is not Gorenstein and
$p_g=1$. If $2P=2Q$, then $p_g\geq2$, and $p_g=2$ if $P\neq Q$. If
$P=Q$ both cases $p_g=2$ and $p_g=3$ are possible.

In the following we give equations for all cases. For double points
$z^2=f(x,y)$ this can be done starting from the topological type of
the branch curve, but in general it is very difficult to find
equations. The easiest approach is to first compute generators of
the local ring. This works for quasi-homogeneous singularities,
where one can use Pinkham's method \cite{pi}. By computing
deformations of positive weight one finds other singularities. As
the graph is star-shaped, there exist quasi-homogeneous
singularities with this graph. Their coordinate ring  is the graded
ring $R=\bigoplus H^0(\sier{E_0}(nQ-\lceil\frac23n\rceil P))$. We
may take the point $Q$ as origin of the group law on the elliptic
curve. The structure of the ring $R$ depends on the order of the
point $P$ on the elliptic curve.

\subsubsection{\boldmath$P=Q\,$, $\;p_g=3\,$.}
In this case the singularity is Gorenstein by a result of S.S.-T.
Yau, see \cite[3.5]{nem}, and all singularities are deformations of
the quasi-homogeneous ones. The weighted homogeneous ring $R$ has
three generators, in degrees 1, 6 and 9. The whole equisingular
stratum is given by
\[
x^2+y^3+a_0yz^{12}+z^{18}+a_{-1}yz^{13}+a_{-2}yz^{14}
  +a_{-3}yz^{15}+a_{-4}yz^{16}\;.
\]

\subsubsection{\boldmath$P=Q\,$, $\;p_g=2\,$.}
A singularity of this type is still a deformation of a
quasi-homogeneous singularity, but now a nonnormal one. The
resolution specialises to the resolution of the quasi-homogeneous
singularity and blowing down the exceptional set of the total space
yields a family specialising to a nonnormal singularity with
$\delta= l(\widetilde{\cal O}/{\cal O})=1$, where $\widetilde{\cal
O}$ is the normalisation of the local ring. An example of such a
singularity is provided by the subring  of
$\C\{x,y,z\}/(x^2+y^3+z^{18})$ generated by $x$, $y$, $xz$, $yz$,
$z^2$ and $z^3$. In computing a deformation of positive weight of
this singularity there is no guarantee that the computation stops,
but with some luck we succeeded. For the deformed singularity two
variables can be eliminated and the result of our computation is a
determinantal singularity in $\C^4$, for which we give the
equisingular stratum. This will contain most singularities of this
type, but maybe not all. In new variables $(x,y,z,w)$ we get the
equations
\[
\Rank\begin{pmatrix}
2w-xz&          4y+x^2&z \\
4xy-z^2-4z^4+b_0x^3z^2+(a_{-1}xz+a_0x^2)(z^2-x^3)&2w+xz&y+x^2
\end{pmatrix} \leq 1  \;.
\]
Checking that this singularity has indeed the required resolution
graph can be done by embedded resolution, for which point blow-ups
suffice. Eliminating $w$, i.e., projecting onto the first three
coordinates, gives the non-isolated hypersurface singularity
\[
z^2(z^2-2x^3-6xy+4z^4-b_0x^3z^2+(a_{-1}xz+a_0x^2)(x^3-z^2))
+(4y+x^2)(y+x^2)^2\;.
\]
The singular locus is given by  $z=y+x^2=0$. Let $\omega$ be the
generator of the dualising sheaf. Then $z\omega$ and $(y+x^2)\omega$
lift to regular differential forms on the normalisation. So indeed
$p_g=2$.

This example is related to examples of Okuma \cite[6.3]{okn}.

\subsubsection{\boldmath$P\neq Q\,$, but $2P=2Q\,$.}
Every such  singularity is a deformation of a
quasi-homogenous one. Generators of the weighted homogeneous ring
have degree 2, 3 and 7. The equisingular stratum is
\[
x^2+z(y^4+a_0y^2z^3+z^6)+a_{-1}yz^{6}+a_{-2}y^2z^{5}
  +a_{-4}y^2z^{6}\;.
\]

\subsubsection{\boldmath$2P\neq 2Q\,$.}
Here there are two moduli, one being the modulus of the elliptic
curve, the other the class  of the divisor $P-Q$ in the Jacobian of
the curve. We start again from the quasi-homogeneous case. The
dimension of the group  $H^0(\sier{E_0}(nQ-\lceil\frac23n\rceil P))$
is $\lfloor\frac n3\rfloor$. We consider the  plane cubic
\begin{equation}\label{cub}
\xi^2\eta+\eta^2\xi+\zeta^3+\la \xi\eta\zeta
\end{equation}
with inflexion point  $Q=(1:0:0)$ and inflexional tangent $\eta=0$.
A section of $H^0(\sier{E_0}(3Q-2 P))$ can be given by the rational
function $\frac l\eta$, where $l$ is the equation of the tangent
line in the point $P$. For $H^0(\sier{E_0}(4Q-3 P))$ we take the
function $\frac q{\eta^2}$, where $q$ is a conic intersecting the
the curve in $Q$ with multiplicity 2 and in $P$ with multiplicity 3.
Continuing in this way one gets generators of the local ring. We
refrain from giving complicated general formulas. The formulas
simplify if $3P=3Q$, but on the other hand the embedding dimension
increases by one. Following the principle of first describing the
most special singularity we take for $P$ the inflexion point
$(0:1:0)$. The function $\frac \xi\eta$ is a section of
$H^0(\sier{E_0}(3Q-2 P))$ and $H^0(\sier{E_0}(4Q-3 P))$. We
homogenise and write the sections as forms of degree $k$. With this
convention we find generators $x=\xi\eta^2$, $y=\xi\eta^3$,
$z=\zeta\xi\eta^3$, $u=\zeta\xi\eta^4$, $v=\zeta^2\xi\eta^4$ and
$w=\zeta^3\xi\eta^5$. We need a generator in degree 9, as $xu=yz$.
There are 10 equations, which can be written in {\em rolling factors
format\/} \cite{re, str}. There are two deformations of positive
weight, which respect this format; in fact, they can be written as
deformation of the equation \eqref{cub}:
\begin{equation}\label{cuba}
\xi^2\eta+\eta^2\xi+\zeta^3+\la \xi\eta\zeta +a_{-1}\xi^2\eta^2
+a_{-3}\xi^2\eta^3\zeta\;.
\end{equation}
The equations consist of the six equations of a scroll:
\[
\Rank\begin{pmatrix}  x&z &v&y\\  z&v& w& u\end{pmatrix} \leq 1
\]
and four additional equations, obtained by multiplying the equation
(\ref{cuba})   with suitable factors. The transition from one
equation to the next involves the replacement of one occurring entry
of the top row of the matrix by the one standing below it. We get
\begin{gather*}
y^3+x^4+xw+\la xyz+a_{-1}yx^3+a_{-3}ux^3\;, \\
uy^2+zx^3+zw+\la yz^2+a_{-1}yx^2z+a_{-3}ux^2z\;, \\
u^2y+vx^3+vw+\la yzv+a_{-1}yx^2v+a_{-3}ux^2v\;, \\
u^3+wx^3+w^2+\la yzw+a_{-1}yx^2w+a_{-3}ux^2w\;.
\end{gather*}
Furthermore there is a deformation of degree 0, corresponding to
moving the point $P$. It does not respect the rolling factors
format; as infinitesimal deformation we change the equation of
degree 9 to $ux-yz+\ep w$, showing again that the embedding
dimension drops for $3P\neq 3Q$. We did not succeed in computing the
higher order terms, which seem to involve power series. Only for the
special value $\la=0$ of the modulus the computation ended after a
finite number of steps. We give the result.
\begin{gather*}
\Rank\begin{pmatrix} z &v&y\\  v& w& u-\ep x^2\end{pmatrix} \leq
1\;,
\\
(u-\ep^4x^2)x-zy+\ep w\;,\\
z^2-vx+\ep(u-\ep^4x^2)y-2\ep^3vx\;,\\
vz-wx+\ep(u-\ep x^2)(u-\ep^4x^2)-2\ep^3 wx\;,\\
y^3+wx+x^4-\ep u^2\;,\\
wz+(u-\ep^4x^2)y^2+zx^3-2\ep^2vyx-z\ep^3x^3\;,\\
wv+(u-\ep x^2)(u-\ep^4x^2)y+vx^3-2\ep^2wyx-v\ep^3x^3\;,\\
w^2+(u-\ep x^2)^2(u-\ep^4x^2)+wx^3-2\ep^2wx(u-\ep x^2)-w\ep^3x^3\;.
\end{gather*}

Similar examples exist with $E_0$ replaced by other minimally
elliptic cycles, such as a chain of rational curves. But if $E_0$ is
a rational cuspidal curve, there are no torsion points, as the group
structure of a cuspidal cubic is that of the additive group
${\mathbb G}_a$, and there is only one possibility for a Gorenstein
singularity, the one with $p_g=3$, giving the singularity $E_{36}$
(in Arnold's notation); a quasihomogeneous  equation is
$x^2+y^3+z^{19}$.

\section{The conjectures}

The conjecture has been made that the geometric genus is
topological, if one assumes  a topological restriction, that the
link of the singularity is a rational homology sphere, and an
analytical one, that the singularity is $\Q$-Gorenstein. In fact,
N\'emethi and Nicolaescu proposed a precise formula:
\[
p_g={\bf sw}(M)-(K^2+r)/8\;,
\]
where ${\bf SW}(M)$ is the Seiberg-Witten invariant of the link $M$,
as before $K$ is the canonical cycle and $r$ the number of
components of the exceptional divisor; for details see \cite{nn}.

For smoothable Gorenstein singularities this conjecture can be
formulated in terms of the signature of the Milnor fibre. We recall
Laufer's formula
\[
\mu =12p_g-b_1(E)+b_2(E)+K^2\;.
\]
Together with $\mu_0+\mu_+=2p_g$ and $\mu_0=b_1(E)$ we find Durfee's
signature formula
\[
\sigma=8p_g-b_2(E)-K^2\;.
\]
So the conjecture is ${\bf SW}(M)=-\sigma/8$.
This conjecture generalises the Casson Invariant Conjecture of
Neumann and Wahl \cite{nwc}, that for an complete intersection with
integral homology sphere link $\sigma/8$ equals the Casson
invariant.

As an outgrowth of  their work Neumann and Wahl even came up with a
way to write down equations from the resolution graph. What is known
about it, is very well described in Wahl's survey \cite{wahl}. In
the case of a rational homology sphere link satisfying certain
conditions  one gets equations together with an action of the finite
group $H=H_1(M,\Z)$. Given a singularity $(X,p)$ with $H$ finite,
the maximal unramified abelian cover of $X\setminus\{p\}$ is a
Galois cover with covering transformation group $H$, which can be
completed with one point. We get a map $(\wt X,p)\to(X,p)$, which is
called the universal abelian cover of $(X,p)$. The recipe for the
equations involve the splice diagram of the singularity, see
\cite{nw,nwu}. Here we recall only the case of a resolution graph of
a quasi-homogeneous singularity, already given in \cite{neu}. There
is one central curve with $k$ arms, which are resolution graphs of
cyclic quotient singularities with a group of order $\al i$,
$i=1,\dots,k$. The associated $k-2$ equations define a
Brieskorn-Pham complete intersection: one has linear equations in
$x_i^{\al i}$ with general enough coefficient matrix. The group $H$
acts diagonally. In general these equations are the building blocks,
which have to be spliced together, leading to splice type equations.
This requires conditions on the graph, which are called the
`semigroup condition' (in order to be able to splice) and the
`congruence condition' (to have an action of $H$). The resulting
singularity is said to be a splice quotient.

\section{Superisolated singularities}
A surface singularity $f\colon (\C^3,0) \to (\C,0)$ is called
superisolated if blowing up the origin once resolves the
singularity. We investigate what this means in terms of equations in
a more general situation.

\begin{propo}\label{super}
Let $f=f_d+f_{d+k}+\cdots$ be the decomposition in  homogeneous
parts. The first blow up has only isolated singularities of
suspension type, of the form $z^k=g_i(x,y)$, where the $g_i(x,y)$
are the local equations of the singularities of the exceptional
curve (in $\P^2$), if and only if the homogeneous polynomial $f_d$
defines a plane curve with isolated singularities, through which the
curve $f_{d+k}=0$ does not pass.
\par
In this situation the resolution graph of the singularity depends
only on $d$, $k$ and the topological type of the singularities
$g_i$.
\end{propo}

Singularities satisfying the conditions of the proposition are said
to be of {\em Yomdin type} \cite{ALM}, as the most important example
is a Yomdin series $f_d+l^{d+k}$, where $l$ is a general linear
form.

\begin{proof}
We compute the blow-up. We may assume that $(0:0:1)$ is a singular
point of the plane curve $f_d$ and compute in the chart with
coordinates $(\xi,\eta,z)$, where $(x,y,z)=(\xi z,\eta z,z)$. Then
the strict transform of $f$ is
\[
\bar f(\xi,\eta,z)=f_d(\xi,\eta,1)+z^k(f_{d+k}(\xi,\eta,1)
+zf_{d+k+1}(\xi,\eta,1) +\dots)\;.
\]
The origin is an isolated singularity of $\bar f$ if and only if
$(0:0:1)$ is an isolated singularity of $f_d$. This singularity is a
$k$-fold suspension singularity if and only if $f_{d+k}(\xi,\eta,1)
+zf_{d+k+1}(\xi,\eta,1) +\dots$ is a unit, so $f_{d+k}(0,0,1)\neq
0$. To determine the embedded resolution it suffices to consider the
singularity $f_d(\xi,\eta,1)+z^k$. The other singular points of
$f_d$ can be treated in a similar way.
\end{proof}

In particular, for a superisolated singularity the first blow up
gives the minimal resolution. It has the curve $C\colon f_d=0$ as
exceptional divisor with self intersection $-d$. In fact, the normal
bundle of $C$ is $\sier C(-1)$, where $H=\sier C (1)$ is the
hyperplane bundle of the embedding of $C$ in $\P^2$. The link $M$ of
the singularity is a rational homology sphere if $C$ is homeomorphic
to $S^2$ (we restrict ourselves to the case of irreducible curves).
This implies that it is locally irreducible, so $C$ is an
irreducible rational curve with only (higher) cusp singularities.
Then $H_1(M,\Z)=\Z/d\Z$, generated by a small loop around the
exceptional curve. The universal abelian cover of such a
superisolated singularity is a cyclic cover of degree $d$. We can
describe it in the following way \cite{ok}. Consider the
$\Q$-divisor $\frac1d C$ on the minimal resolution $\pi \colon \wt
X\to X$ of the singularity. There exists a divisor $D$ on $\wt X$
such that $C$ is linearly equivalent to $dD$. Then $\pi_*D$ is a
$\Q$-Cartier divisor of index $d$ on $X$. The local ring of the
covering $Y$, as $\sier X$-module, is $\bigoplus_{j=0}^{d-1} \sier
X(-j\pi_*D)=\bigoplus_{j=0}^{d-1} \pi_*\sier {\wt X}(-jD)$. The
sheaf $\bigoplus_{j=0}^{d-1} \sier {\wt X}(-jD)$ defines a $d$-fold
cover of $\wt X$, branched along $C$, which is a partial resolution
$\overline Y$ of $Y$. It has $C$ as exceptional divisor, this time
with self intersection $-1$, and singularities at the singular
points of $C$, of suspension type $z^d=g_i(x,y)$. It is now easy to
find the resolution graph of the universal abelian cover, as already
noted in \cite[3.2]{lmn}. But is is very hard to find $D$ and to
determine the ring structure explicitly. We therefore compute on the
curve $C$.

The exceptional curve $C$ on the partial resolution $\overline Y$
has as normal bundle the dual of a $d$th root of $H=\sier C (1)$.
For a rational cuspidal curve this $d$th root is unique. In general,
for a singular curve $C$ with normalisation $n\colon \wt C \to C$,
there is an exact sequence \cite[Exercise II.6.8]{hag}
\[
0\lra \textstyle\bigoplus_{P\in C} \wt{\cal O}_P^*/\sier P ^* \lra
\Pic C \lra \Pic \wt C \lra 0\;,
\]
where $\wt{\cal O}_P$ is the normalisation of the local ring of $P$,
and the star denotes the group of units in a ring. For a singularity
with one branch $\wt{\cal O}_P^*/\sier P ^*$ is a unipotent
algebraic group. With $\Pic \wt C=\Z$ for a rational curve this
gives the uniqueness.

\section{Degree four}
In this section we construct examples of universal abelian covers,
which are not complete intersections. To this end we consider
superisolated singularities with lowest degree part of degree four.
The singularities on  cuspidal rational quartic curves can be $E_6$,
$A_6$, $A_4+A_2$ or $3A_2$. This list is easily obtained, as the sum
of the Milnor numbers has to be 6. Equations can be found in
\cite{sal}.

\subsection{Rational curves with an \boldmath$E_6$-singularity}
There are two projectively inequivalent quartic curves with an
$E_6$-singularity. One admits a $\C^*$-action, and can be given by
the equation $x^4-y^3z=0$. The other type can be written as
$(x^2-y^2)^2-y^3z=0$. In the first case the line $z=0$ is a
hyperflex, while it is an ordinary bitangent in the second case.

The first type gives rise to a series of quasi-homogeneous isolated
surface singularities
\[
f_k\colon \; x^4-y^3z-z^{4+k}\;.
\]
In particular, for $k=1$ we have a superisolated singularity, whose
minimal resolution has the given rational curve as exceptional
divisor with self-intersection $-4$. The minimal good resolution has
dual graph:
\[
\unitlength=30pt
\def\ci{\circle*{0.26}}
\def\vi{\ci}
\def\mbt#1{\makebox(0,0)[t]{$\scriptstyle #1$}}
\begin{picture}(3,1.5)(-1.26,-.25)
\put(-1,0){\line(1,0){.87}}  \put(2,0){\line(-1,0){1.87}}
\put(-1,0){\vi} \put(-1,-.2){\mbt{-16}} \put(1,0){\ci}
\put(0,1){\line(0,-1){.87}} \put(0,1){\vi}
\put(0,1.2){\makebox(0,0)[b]{$\scriptstyle -4$}} \put(2,0){\ci}
\put(0,0){\circle{0.26}} \put(0,-.2){\mbt{-1}}
\end{picture}
\]
A splice type equation for this graph is $\xi^4-\eta^3-\zeta^{16}$.
The action of $H_1(M,\Z)$ is $\frac14(3,4,1)$. This means that we
have a diagonal action of $\Z_4$ with generator
$(\xi,\eta,\zeta)\mapsto (-i\xi,\eta,i\zeta)$. Invariants are
$x=\xi\zeta$, $y=\eta$ and $z=\zeta^4$ (note that
$\xi^4=\eta^3+\zeta^{16}$) and the quotient is exactly the
quasi-homogeneous superisolated singularity. In fact, the same group
acts on the singularities $\xi^4-\eta^3-\zeta^{12+4k}$ with quotient
exactly $f_k$, but in general this will not be the universal abelian
cover. The link need not even be a rational homology sphere (an
example in case is $k=12$).

For $k=2$ the graph is:
\[
\unitlength=30pt
\def\ci{\circle*{0.26}}
\def\vi{\ci}
\def\mbt#1{\makebox(0,0)[t]{$\scriptstyle #1$}}
\begin{picture}(4,2.25)(-2.26,-1)
\put(-2,0){\line(1,0){4}} \put(-1,0){\ci} \put(-2,0){\ci}
\put(1,0){\ci} \put(0,1){\line(0,-1){2}} \put(0,1){\vi}
 \put(0,-1){\ci}
\put(0,1.2){\makebox(0,0)[b]{$\scriptstyle -10$}} \put(2,0){\ci}
\put(0,0){\ci}
\end{picture}
\]
The quasi-homogeneous singularity $f_2$ has  a universal abelian
cover with splice type equations of the form
\begin{align*}
  z_1^{10} \;&=\; az_2^{3}+bz_3^{3} \\
  z_4^{2} \;&=\; cz_2^{3}+dz_3^{3}
\end{align*}
with diagonal  action of a group of order 12, generated by
$[-i,-1,-1,i]$ and $[-1,\ep^2,\ep,-1]$, where $\ep$ is a primitive
third root of unity.

The second type of quartic curve with an $E_6$-singularity, the one
with ordinary bitangent,  gives a series of singularities with the
same topology as the series $f_k$ (see Prop.~\ref{super}):
\[
g_k\colon \; (x^2-y^2)^2-y^3z-z^{4+k}\;.
\]
But this time the fourfold cover is not a complete intersection.
There is no easy recipe for the equations. We  look at the  cone
over the singular quartic curve $C$. As a quartic curve it is
canonically embedded. Let $K$ be the canonical line bundle on $C$.
There is a unique line bundle $L$ of degree one on $C$ with $4L=K$.
We determine the ring $\bigoplus H^0(C,nL)$. We have $H^0(C,L)=0$:
an effective divisor of degree $1$ would give an hyperflex. But
$h^0(C,2L)=1$, the effective divisor of degree $2$ being the points
of tangency of the bitangent. By Riemann-Roch $h^0(C,3L)=1$. So we
have to find the unique effective divisor $D$ of degree 3, for which
there is a cubic curve with fourfold contact. The cubic lies in the
system of contact cubics, which we obtain (as in the classical case
of smooth quartics \cite[\S 17]{kl}) by writing the equation of the
curve in the form $Q^2-Tz$, where $Q=x^2-y^2-z(ax+by+cz)$ and
$T=y^3-2(x^2-y^2)(ax+by+cz)+z(ax+by+cz)^2$. As the curve $C$ is
invariant under the involution $x\mapsto-x$, the cubic (being
unique) is also invariant, from which we conclude that $a=0$. Using
the parametrisation $(x,y,z)=(st^3,t^4,(s^2-t^2)^2)$ we get
$T=\left[t^6-(s^2-t^2)(bt^4+c(s^2-t^2)^2)\right]^2$. Writing out the
condition that the term in brackets is the square of a cubic form
and solving gives $b=3$, $c=-4$. So
\begin{align*}
  Q \;&=\; x^2-y^2-z(3y-4z) \;,\\
  T \;&=\; y^3-2(x^2-y^2)(3y-4z)+z(3y-4z)^2 \;.
\end{align*}
The three points in $D$ are $(0:1:1)$ and $(\pm2\sqrt6:4:1)$ or in
terms of the parametrisation $(0:1)$ and $(\pm\sqrt3:\sqrt2)$. The
curve $C$ is the boomerang shaped curve (with the $E_6$ singularity
at the tip) in the next figure. To make the bitangent visible the
coordinate transformation $z\mapsto z-\frac12 y$ was applied. The
affine chart $z=1$ in the new coordinates is shown.  The conic $Q$
intersects $C$ in the points of tangency of the bitangent and is
tangent in the three points of contact with the contact cubic $T$,
whose real locus consists of a very small, hardly visible oval, and
an odd branch with one point at infinity (an inflexion point); the
picture shows two arcs of this branch.
\begin{center}
\psfrag{T}{$T$} \psfrag{C}{$C$} \psfrag{Q}{$Q$}
\includegraphics[width=.9\linewidth]{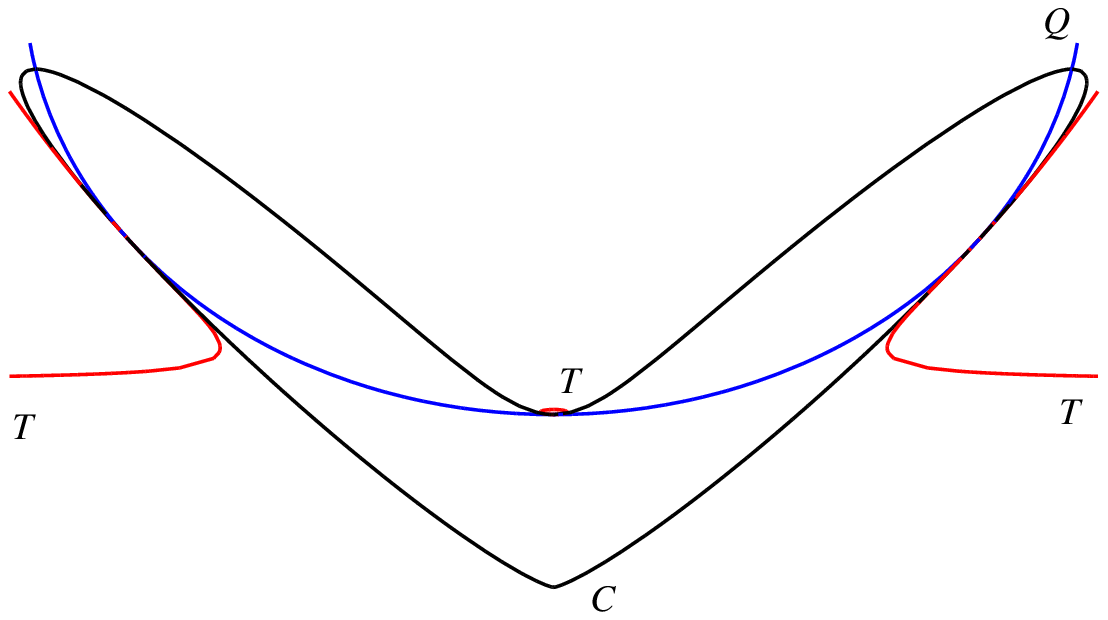}
\end{center}
The picture was made with Richard Morris' program {\sc
SingSurf} \cite{SS}.

As $4D=3K$ we have indeed that $L=K-D$ satisfies $4L=K$. We can
describe $H^0(C,n(K-D))$ with polynomials of degree $n$, passing $n$
times through $D$. In degree 2 we have $Q$, in degree 3 $T$, in
degree four $Tx$, $Ty$ and $Tz=Q^2$. One section in degree 5 is
$TQ$, the others are of the form $TQ'$ with $Q'$ a quadric through
the three points in $D$. Alternatively we can write forms of degree
$n$ in $s$ and $t$. The sections described so far generate the ring.
As minimal set we choose the following generators:
\[
\begin{array}{c@{\colon\quad}l@{\quad=\quad}l}
  \zeta & Q             & s^2-t^2 \\
  u     & T             & s(2s^2-3t^2) \\
  x     & Tx            & st^3 \\
  y     & Ty            & t^4 \\
  v     & T(y-4z)x      & -(2s^2-t^2)t^3 \\
  w     & T(y-4z)(y-z)  & s(2s^2-t^2)(s^2-2t^2) \\
\end{array}
\]
One finds 9 equations between the 6 variables. They can be given in
{\em rolling factors format\/} \cite{str}.  The first 6 equations
are the $2\times 2$ minors of the matrix
\[
\addtocounter{equation}{1}\label{eza}
\begin{pmatrix} u & y-\zeta^2 & x & w \\
y-4\zeta^2 & w & v & u^2-8\zeta^3+6y\zeta\end{pmatrix}\;,
\tag{\theequation.a}
\]
while the remaining 3 equations come from $u^2\zeta=Q(x,y,\zeta^2)$
by rolling factors. Specifically,
\begin{gather*}
  u^2\zeta-x^2+(y-\zeta^2)(y+4\zeta^2)\;,\\
  (y-4\zeta^2)u\zeta-vx+w(y+4\zeta^2)\;, \label{ezb}\tag{\theequation.b}\\
  (y-4\zeta^2)^2\zeta-v^2+(u^2-8\zeta^3+6y\zeta)(y+4\zeta^2)\;.
\end{gather*}

By perturbing these equations we get again a series of
singularities. We use a rolling factors deformation \cite{str}: we
deform the matrix and use the deformed matrix for rolling factors.
As there is a recipe for the relations, basically in terms of the
matrix, this yields a (flat) deformation of the singularity. The
easiest way to perturb is to change the last entry of the matrix
(\ref{eza}) into
\begin{equation}\label{ezphi}
    u^2-8\zeta^3+6y\zeta+\vp_k
\end{equation}
and correspondingly the term in the third additional equation of
(\ref{ezb}). We take $\vp_{2k-1}=u\zeta^{1+k}$ and
$\vp_{2k}=\zeta^{3+k}$, because the singular locus of the quasi-cone
$\oplus H^0(C,nL)$, $4L=K$, is given by $u^2=4\zeta^3$, $w=u\zeta$.
As the ring $\oplus H^0(C,nL)$ is Gorenstein by construction, all
the isolated singularities we obtain are Gorenstein. A weighted
blow-up (with the weights of the non-isolated singularity) gives a
partial resolution with the curve $C$ as exceptional divisor and one
singularity of suspension type $z^k=x^3+y^4$.

By choosing an equivariant deformation for the $\Z_4$-action we
ensure that it descends to the quotient. As quotient of the deformed
singularities with $\vp_{4k}$ we obtain the Yomdin type series
\[
\widetilde g_k\colon \;
(x^2-y^2)^2-y^3z+(x^2-y^2-z(3y-4z))z^{2+k}\;.
\]
For $k=1$ we have indeed a superisolated singularity. One computes
that the Tyurina number of this singularity is 28, whereas
$\tau(g_1)=30$. So we did not compute the universal abelian covering
of $g_1$, but of an analytically inequivalent singularity with the
same degree 4 part. The general philosophy of superisolated
singularities dictates that the  covering of $g_1$ is a deformation
of the same quasi-cone $\oplus H^0(C,nL)$. As it is not easy to give
a direct formula, we start again with a rolling factors deformation,
this time changing the $(2,1)$-entry of the matrix (\ref{eza}) into
$y-4\zeta^2+16\zeta^4+20y\zeta^2+21y^2+216\zeta^6$  and
correspondingly the additional equations (\ref{ezb}). We get  a
quotient of the form $(x^2-y^2)^2-\upsilon_1y^3z-\upsilon_2z^5$,
where $\upsilon_1$ and $\upsilon_2$ are units, and an easy
coordinate change brings us to $g_1$. The precise formulas for the
units are rather complicated.

The computations so far can be used to give other  examples. A
related series is obtained by first taking the double cover. For a
plane quartic with bitangent $z=0$ with equation of the form
$Q^2-Tz$ the ring $\oplus H^0(C,n\Theta)$, where $\Theta$ is the
line bundle belonging to the bitangent, has four generators, $\zeta$
in degree 1, with $\zeta^2=z$, in degree 2 generators $x$ and $y$
and finally in degree 3 a generator $w=\sqrt T$. Equations are
\begin{align*}
  \zeta w \;&=\; Q(x,y,\zeta^2)\;, \\
  w^2 \;&=\;T(x,y,\zeta^2)\;.
\end{align*}
We obtain a series of singularities from our quartics with $E_6$ by
changing the second equation into $w^2 = T(x,y,\zeta^2)+\psi_k$, the
easiest being $\psi_k=\zeta^{6+k}$. We look at the first element of
the series.

The family
\begin{align*}
  \zeta w \;&=\; x^2-ty^2\;, \\
  w^2 \;&=\; y^3+\zeta^7\;
\end{align*}
is equisingular. Its minimal resolution has a  rational curve with
an $E_6$ of self-intersection $-2$ as exceptional divisor. The
universal abelian cover is  a double cover. For $t=0$ it is the
hypersurface singularity $\xi^4-\eta^3-\zeta^{14}$, while for
$t\neq0$ it has embedding dimension 6, being an appropriate
deformation of our singular quasi-cone with equations (\ref{eza})
and (\ref{ezb}); taking $\vp_2=\zeta^4$ in the perturbation
(\ref{ezphi}) yields a quotient with a more complicated formula.

Finally, by looking at the lowest element of the series  for the
fourfold cover we get an example of a Gorenstein singularity with
integral homology sphere link, which is not a complete intersection.

\begin{propo} The Gorenstein singularity obtained by taking
$\vp_1=u\zeta^2$  in (\ref{ezphi}), i.e., with equations
\begin{gather*}
 \Rank\begin{pmatrix} u & y-\zeta^2 & x & w \\ y-4\zeta^2 & w & v &
   u^2-8\zeta^3+6y\zeta+u\zeta^2\end{pmatrix}\leq 1\;,\\[\smallskipamount]
   u^2\zeta-x^2+(y-\zeta^2)(y+4\zeta^2)\;,\\
  (y-4\zeta^2)u\zeta-vx+w(y+4\zeta^2)\;, \\
  (y-4\zeta^2)^2\zeta-v^2+(u^2-8\zeta^3+6y\zeta+u\zeta^2)(y+4\zeta^2)\;,
\end{gather*}
has integral homology sphere link. Its minimal resolution has the
rational curve  with an $E_6$-singularity as excep\-tional divisor
with self-intersection $-1$. The resolution graph is the same as for
the hypersurface singularity $\xi^4-\eta^3-\zeta^{13}$\/{\rm:}
\[
\unitlength=30pt
\def\ci{\circle*{0.26}}
\def\vi{\ci}
\def\mbt#1{\makebox(0,0)[t]{$\scriptstyle #1$}}
\begin{picture}(3,1.5)(-1.26,-.25)
\put(-1,0){\line(1,0){.87}}  \put(2,0){\line(-1,0){1.87}}
\put(-1,0){\vi} \put(-1,-.2){\mbt{-13}} \put(1,0){\ci}
\put(0,1){\line(0,-1){.87}} \put(0,1){\vi}
\put(0,1.2){\makebox(0,0)[b]{$\scriptstyle -4$}} \put(2,0){\ci}
\put(0,0){\circle{0.26}} \put(0,-.2){\mbt{-1}}
\end{picture}
\]
\end{propo}

For the hypersurface singularity $p_g=8$. There are several ways to
compute this. One can determine the number of spectrum numbers less
or equal to one, an easy computation for a Brieskorn polynomial ---
a lazy method is to use the computer algebra system {\sc Singular}
\cite{GPS} to find the spectrum. Or one can use the Laufer formula
$\mu=12p_g+K^2+b_2(E)-b_1(E)$ on the minimal resolution. As
explained at the beginning of the next section, the geometric genus
can also be computed as $\sum _{n=0}^4h^0(C,nL)$ on the singular
curve. For the non hypersurface singularity one finds in this way
$p_g=6$.

\subsection{The three-cuspidal quartic}
\label{cusp}

A quartic curve with three $A_2$-singularities  has a bitangent, and
equations for the fourfold cover can be found in the same way as in
the case of a curve with an $E_6$ singularity and  an ordinary
bitangent.

Let $\sigma_1$, $\sigma_2$ and $\sigma_3$ be the elementary
symmetric functions of $x$, $y$ and $z$. The 3-cuspidal curve has
equation $\sigma_2^2-4\sigma_1\sigma_3=
x^2y^2+x^2z^2+y^2z^2-2x^2yz-2xy^2z-2xyz^2$. The bitangent is
$\sigma_1=x+y+z$. We have again $L=K-D$, with $D$ a divisor of
degree $3$, consisting of the points $(1:4:4)$, $(4:1:4)$ and
$(4:4:1)$. We write the curve in the form $Q^2-4T\sigma_1$, where
the cubic form $T$ cuts out $4D$. One computes
\[
\begin{array}{c@{\quad=\quad}c@{\quad=\quad}c}
  Q   &  27\sigma_2-8\sigma_1^2
    & -8(x^2+y^2+z^2)+11(xy+xz+yz)\;,\qquad\qquad\\[\smallskipamount]
  T   &   729\sigma_3-108\sigma_1\sigma_2+16\sigma_1^3\\
   \multispan{3} $\hfill
   16(x^3+y^3+z^3)-60(x^2y+xy^2+x^2z+y^2z+xz^2+yz^2)+501xyz\;.\hfill$
\end{array}
\]

The quadric $Q$ is a section of $2L$ and the cubic $T$ a section of
$3L$. A basis of sections of $4L=4K-4D\equiv K$ consists of $Tx$,
$Ty$ and $Tz$; due to the relation $Q^2=4T\sigma_1$ we need only two
of these as generators of the ring. The necessary choice breaks the
symmetry. In degree five we have sections of the form $TQ'$, where
$Q'$ is a quadric passing though the three points of the divisor
$D$. We choose as generators two forms with reducible quadrics,
consisting of the line through two of the points and a line though
the third point. The equations are again in rolling factors format,
with the transition from the first to the second row being
multiplication with the equation of the chosen line. A suitable
choice for the other line in a generator of degree 5 is the tangent
in the third point. Therefore we take as generators $\frac12Q$, $T$,
$Ty$, $Tz$, $T(y-z)(5x-4y-4z)$ and $T(8x-y-z)(5x-4y-4z)$, which
correspond to coordinates $(s,u,y,z,v,w)$.

We get a series of isolated singularities by a rolling factors
deformation. With the given generators we form the matrix
\[\addtocounter{equation}{1}\label{dca}
\begin{pmatrix} u & y-z & 8s^2-9y-9z & w \\ 5s^2-9y-9z & v & w &
u^2-46s^3+54ys+54zs+\vp_k\end{pmatrix}\;,\tag{\theequation.a}
\]
while the last 3 equations come from $u^2s=Q(s^2-y-z,y,z)$ by
rolling factors. Specifically,
\begin{gather*}
  8u^2s+(8s^2-9y-9z)(4s^2-9y-9z)+27(y-z)^2\;,\\
  8(5s^2-9y-9z)us+w(4s^2-9y-9z)+27(y-z)v\;, \tag{\theequation.b}\\
  8(5s^2-9y-9z)^2s+(u^2-46s^3+54ys+54zs+\vp_k)(4s^2-9y-9z)+27v^2\;.
\end{gather*}
We take $\vp_{2k-1}=us^{1+k}$ and $\vp_{2k}=s^{3+k}$, while
$\vp_\infty=0$ gives the undeformed non-isolated singularity. A
partial resolution exists with the curve $C$ as exceptional divisor
and three singularities of suspension type $z^k=x^2+y^3$.

With $k=4$ we get the universal abelian cover of a superisolated
singularity. One gets a nicer quotient by throwing in some factors
9: the above equations with $\vp_4=729s^5$ give the universal
abelian cover of the superisolated singularity
$\sigma_2^2-4\sigma_1\sigma_3-2Q\sigma_1^3$, which has a 3-cuspidal
rational curve with self-intersection $-4$ as exceptional divisor on
the minimal resolution.

\begin{propo}\label{cuspprop}
The Gorenstein singularity in $(\C^6,0)$ with  equations
(\ref{dca},b) with $\vp_1=us^2$ has integral homology sphere link.
The exceptional divisor on the minimal resolution is the
three-cuspidal rational curve with self-intersection $-1$. The
resolution graph and splice diagram are as follows. The semigroup
condition is not satisfied.
\[
\unitlength=30pt
\def\ci{\circle*{0.26}}\def\cir{\circle{0.26}}
\def\mb#1#2{\makebox(0,0)[#1]{$\scriptstyle #2$}}
\begin{picture}(5,3.5)(-2.5,-2.5)
\put(0,0){\ci} \put(0,.2){\mb b{-19}}
\put(-1.5,0){\cir} \put(-1.5,-.2){\mb t{-1}}
\put(-2.5,-1){\line(1,1){.9}}
\put(-2.5,-1){\ci} \put(-2.5,-1.2){\mb t{-3}}
\put(-2.5,1){\line(1,-1){.9}} \put(-2.5,1){\ci}
\put(0,0){\line(-1,0){1.37}}
\put(1.5,0){\cir} \put(1.5,-.2){\mb t{-1}}
\put(2.5,-1){\line(-1,1){.9}}
\put(2.5,-1){\ci} \put(2.5,-1.2){\mb t{-3}}
\put(2.5,1){\line(-1,-1){.9}} \put(2.5,1){\ci}
\put(0,0){\line(1,0){1.37}}
\put(0,-1.5){\cir} \put(0.2,-1.5){\mb l{-1}}
\put(-1,-2.5){\line(1,1){.9}}
\put(-1,-2.5){\ci} \put(-1.2,-2.5){\mb r{-3}}
\put(1,-2.5){\line(-1,1){.9}} \put(1,-2.5){\ci}
\put(0,0){\line(0,-1){1.37}}
\end{picture}
\qquad
\def\ci{\circle{0.2}}
\begin{picture}(6.5,3.5)(-3.25,-3)
\put(0,0){\ci}
\put(.1,0){\line(1,0){1.3}}  \put(1.2,.1){\mb b7} \put(.3,.1){\mb b1}
\put(1.5,0){\ci}
\put(1.6,.04){\line(3,1){1.3}} \put(1.8,.22){\mb b2}
\put(1.6,-.04){\line(3,-1){1.3}} \put(1.8,-.22){\mb t3}
\put(3,.505){\ci}
\put(3,-.505){\ci}
\put(-.1,0){\line(-1,0){1.3}}  \put(-1.2,.1){\mb b7} \put(-.3,.1){\mb b1}
\put(-1.5,0){\ci}
\put(-1.6,.04){\line(-3,1){1.3}} \put(-1.8,.22){\mb b2}
\put(-1.6,-.04){\line(-3,-1){1.3}} \put(-1.8,-.22){\mb t3}
\put(-3,.505){\ci}
\put(-3,-.505){\ci}
\put(0,-.1){\line(0,-1){1.3}}  \put(.1,-1.2){\mb l7} \put(.1,-.3){\mb l1}
\put(0,-1.5){\ci}
\put(.04,-1.6){\line(1,-3){.435}} \put(.22,-1.8){\mb l2}
\put(-.04,-1.6){\line(-1,-3){.435}} \put(-.22,-1.8){\mb r3}
\put(.505,-3){\ci}
\put(-.505,-3){\ci}
\end{picture}
\]
\end{propo}
Indeed, the semigroup condition requires that a 1 adjacent to the
central node is in the semigroup generated by 2 and 3, which is
impossible.

This singularity is in fact a universal abelian cover of a
hypersurface singularity, which will be considered in section
\ref{verder}.

\subsection{The case \boldmath$A_6$}
This is one of the examples made more explicit by  Luengo-Velasco,
Melle-Her\-n\'an\-dez and  N\'emethi \cite[4.5]{lmn}. A quartic
curve with $A_6$ is unique up to projective equivalence. Its
equation is $(zy-x^2)^2-xy^3=0$. The associated superisolated
singularity is different from the splice quotient associated to its
resolution graph. In  that case the corresponding curve of
arithmetic genus 3 is the weighted complete intersection $yz=x^2$,
$t^2=xy^3$. This shows that the curve is hyperelliptic: the
canonical linear system is not very ample. Seven Weierstra\ss\
points are concentrated in the singular point, while
$P\colon(x:y:z)=(0:1:0)$ is an ordinary Weierstra\ss\ point. One has
$4P=K$. The ring $\oplus H^0(C,nP)$ is simply given by
$\{u^7=w^2\}\subset \C^3$.

The case  $A_4+A_2$ is similar, in that the splice quotient comes
from an hyperelliptic curve, while the plane quartic is unique with
equation $(xy-z^2)^2-yx^3=0$.

The plane quartic with an $A_6$ is a deformation of the
hyperelliptic curve. As $2P$ is the $g_2^1$, which is an even
theta-characteristic, the plane quartic with an $A_6$ has a unique
ineffective theta-characteristic $\Theta$, which is a line bundle.
Then it is well known, going back to Hesse,
that one can write the equation of the curve as linear
symmetric determinant. It is possible to compute the matrix, but the
easiest approach is to search for it in  Wall's classification of
nets of quadrics \cite{wall}. After some elementary operations we
find the matrix
\[
M=\begin{pmatrix}
y  & x   &  z  &   0  \\
x  & 4z  &  0  &  -y  \\
z  & 0   &  x  &   x  \\
0  & -y  &  x  &   0\end{pmatrix}  \;.
\]
We get a series of singularities by changing the last entry of the
matrix. In particular, replacing the $0$ on the diagonal by $z^2$
leads to the superisolated singularity
\[
(zy-x^2)^2-xy^3+z^2(4yzx-x^3-4z^3)\;.
\]

The matrix $M$ not only gives the curve $C$, as its determinant, but
also the embedding by $K+\Theta=3\Theta$ in $\P^3$: $M$ is the
matrix of a net of quadrics in $\P^3$ and the curve is the Steiner
curve of this net, the locus of vertices of singular quadrics. The
matrix also determines the equations (and in fact the whole
resolution of the ideal) of the ring $\oplus H^0(C,n\Theta)$. This
ring has generators $x$, $y$ and $z$ in degree $2$ and four
generators $v_0,\dots v_3$ in degree $3$. There are 14 equations,
which can be written succinctly as matrix equations
\[
Mv=0\;,\qquad vv^t=\textstyle \bigwedge^3M\;,
\]
where $v$ is the column vector $(v_0,v_1,v_2,v_3)^t$. The same
equations for the perturbed matrix define a double cover of the
superisolated singularity. It has the same topology as the
hypersurface singularity $v^{10}-vu^7+t^2=0$, but different $p_g$.
The link is a rational homology sphere. The resolution graph is
\[
\unitlength=30pt
\def\ci{\circle*{0.26}}
\def\vi{\ci}
\def\mbt#1{\makebox(0,0)[t]{$\scriptstyle #1$}}
\begin{picture}(6,1.5)(-3.26,-.25)
\put(-3,0){\line(1,0){2.87}}  \put(3,0){\line(-1,0){2.87}}
\put(-1,0){\vi} \put(-2,0){\ci} \put(-3,0){\ci}
\put(-1,-.2){\mbt{-3}} \put(1,0){\vi}
\put(0,1){\line(0,-1){.87}} \put(1,-.2){\mbt{-3}} \put(0,1){\vi}
\put(0,1.2){\makebox(0,0)[b]{$\scriptstyle -9$}} \put(2,0){\ci}
\put(3,0){\ci} \put(0,0){\circle{0.26}} \put(0,-.2){\mbt{-1}}
\end{picture}
\]

For the fourfold cover we have to find the  fourth root $L$ of $K$.
We look again for an effective divisor $D$ of degree $3$ such that
$4D=3K$ and $L=K-D$. We start by computing the sections of
$3\Theta=6L$ in terms of the parametrisation
$(x:y:z)=(s^2t^2:t^4:s^4+st^3)$. Then we determine which section is
a perfect square. We find the divisor  given by $4s^3+3t^3$. We now
use the following general procedure to compute the ring $\bigoplus
H^0(C,nL)$ in terms of the parametrisation. Sections of $5L$ are
found by finding the sections of $2K=8L$, which are divisible by
$4s^3+3t^3$. The sections of $7L$ are found from the sections of
$10L$. Knowing the sections in sufficiently many degrees we can
figure out the generators of the ring. There are 12 generators:
\begin{equation}\label{azgen}
\begin{array}{@{\mbox{degree }}c@{\colon\qquad}l}
    3 & 4s^3+3t^3, \\
    4 & s^2t^2,\quad t^4,\quad s^4+st^3, \\
    5 & st^4,\quad 4s^3t^2+t^5,\quad 4s^5+5s^2t^3 ,\\
    6 & t^6,\quad s^2t^4,\quad 2s^4t^2+st^5, \\
    7 & t^7,\quad st^6.
 \end{array}
\end{equation}
A computation shows that there are 54 equations between the
generators. To describe their structure it is better to view the
ring $R=\bigoplus H^0(C,nL)$ as module over $R_0=\bigoplus
H^0(C,nK)=\C[x,y,z]/((zy-x^2)^2-xy^3)$. We have $R=R_0\oplus R_1
\oplus R_2 \oplus R_3$ with $R_i=\bigoplus_n H^0(C,(4n+i)L)$. We
first look at $R_3$. As the section $(4s^3+3t^3)^2$ of $6L$ is a
linear combination of the $v_i$, we can arrange that its square is a
principal minor of a matrix defining the curve $C$. We change the
matrix $M$ into
\begin{equation}\label{azmat}
M=\frac12\begin{pmatrix} 0  & y   &  -x  &   0  \\ y  & -32z  &  -9y
& 8x \\ -x  & -9y   &  10x  &   8z  \\ 0  & 8x &  8z  &   -8y
\end{pmatrix}  \;.
\end{equation}
The first principal minor is
\[
W=-80x^3+81y^3+176xyz+256z^3\;.
\]
This form cuts out $4D$ on $C$. We introduce a dummy variable
$\sigma = 1/\sqrt[4] W$ of weight $-3$ (the variables $x$, $y$ and
$z$ having weight 4), so satisfying $\sigma^4W=1$. Then $w=\sigma^3
W$ has indeed divisor $D$. We compute the ideal of $3D$ in the
homogeneous coordinate ring of $C$; it is generated by $W$ and
\begin{gather*}
   -28x^3y+27y^4+44xy^2z-64x^2z^2+64yz^3\;, \\
   x^2y^2-16x^3z+27y^3z+80xyz^2\;.
\end{gather*}
We obtain variables $r_1$ and $r_2$ of degree 7, having the same
divisors as the generators of degree 7 in (\ref{azgen}), by
multiplying these expressions with $\sigma^3$. The generators of
$R_2$ are $(w^2,v_1,v_2,v_3)$. We obtain them from the matrix M by
$v_i= \sigma^2 (\bigwedge^3M)_{0i}$; the matrix is chosen in such a
way that these sections correspond exactly to the generators in
degree 6, given in (\ref{azgen}) in terms of the parametrisation.
For $R_1$ we look at sections of $5(K-D)\equiv 2K-D$. The linear
system of quadrics through  the three points of $D$ can be computed
as the radical of the ideal generated by the four minors of $M$ of
the first row. We find the ideal
\[
(3y^2+16xz, xy-16z^2, x^2+3yz)\;.
\]
Multiplying these generators by $\sigma$ yields the variables $u_1$,
$u_2$ and $u_3$. The dummy variable $\sigma$ gives us the ring
structure. By eliminating $\sigma$ we find the equations.

The matrix $M$ plays again an important role. Let
$v=(w^2,v_1,v_2,v_3)^t$. Then
\[\addtocounter{equation}{1}\label{aza}
Mv=0\;,\qquad vv^t=\textstyle\bigwedge^3M\;,\tag{\theequation.a}
\]
but some equations are consequences of others. We get also equations
from the syzygies of the ideal $(3y^2+16xz, xy-16z^2, x^2+3yz)$
giving $(u_1,u_2,u_3)$ and from those of the ideal leading to
$(r_1,r_2,w)$. By a suitable choice of generators the same matrix
can be used for both ideals (or rather the matrix and its
transpose). Let $r$ be the vector $(r_1,r_2,w)^t$ and
$u=(u_1,u_2,u_3)^t$ and consider the matrix
\[
N=\frac 12 \begin{pmatrix} -x & 4z & -xy \\ 3y & 4x& -y^2 \\
                    16z & -4y& 4x^2-4yz\end{pmatrix}\;.
\]
Then
\[\label{azb}
u^tN=0\;, \qquad Nr=0\;,\qquad ur^t+\textstyle\bigwedge^2N=0\;.
\tag{\theequation.b}\]
The last expression includes for example the
equation $u_1w=3y^2+16xz$, which obviously holds by our definition
of the generators. The remaining equations concern the rewriting of
monomials and they are in a certain sense a consequence of the given
ones. E.g., the equation $w^3-27yu_1+16zu_2+80xu_3$ can be obtained
because we know how to express $w^4$ and $u_iw$ in terms of $x$, $y$
and $z$. But basically the equation boils down to expressing the
cubic $W=-80x^3+81y^3+176xyz+256z^3$, passing through $D$, in the
generators of the ideal of $D$. The space defined by the matrix
equations (\ref{aza},b) alone has other components, but they all lie
in the hyperplane $\{w=0\}$. We can find the ideal of the
singularity by saturating with respect to the variable $w$.

We obtain a superisolated singularity as determinant of a
perturbation of the matrix $M$: we change the upper left entry of
(\ref{azmat}) into $8z^2$. To find its universal abelian cover we
note that the intersection of the singularity with the cubic cone
$W=0$ still consists of three lines, counted with multiplicity four.
We can therefore compute with the same forms, but now over the local
ring of the singularity. In this instance we succeed in determining
generators of the ring. That we now have a deformation of the
non-isolated singularity follows because we only change the ring
structure on the same underlying $\C$-module. Again we find
equations by eliminating the dummy variable $\sigma$. As a
consequence of our set-up the perturbed matrix $M$ enters, but also
$N$ changes, as its determinant is a multiple of the defining
equation of the hypersurface singularity.

\begin{propo}
The superisolated singularity
\[
(yz-x^2)^2-xy^3+(-80x^3+81y^3+176xyz+256z^3)z^2\;,
\]
with a rational curve with an $A_6$-singularity as exceptional
divisor, has as universal abelian cover a singularity of embedding
dimension 12. 

The universal abelian cover can be given by the following matrix
equations and 30 additional equations:
\[
Mv=0\;,\qquad vv^t=\textstyle\bigwedge^3M\;,\qquad u^tN=0\;, \qquad
Nr=0\;,\qquad ur^t+\bigwedge^2N=0\;.
\]
Here $v=(w^2,v_1,v_2,v_3)^t$, $r=(r_1,r_2,w)^t$, $u=(u_1,u_2,u_3)^t$
are vectors involving nine of the variables and the matrices $M$ and
$N$ depend on $x$, $y$ and $z$:
\[
M=\frac12\begin{pmatrix}8z^2  & y   &  -x  &   0  \\ y  & -32z  &  -9y
& 8x \\ -x  & -9y   &  10x  &   8z  \\ 0  & 8x &  8z  &   -8y\end{pmatrix}
\]
and
\[
N=\frac 12 \begin{pmatrix} -x & 4z & -xy +108yz^2\\ 3y & 4x& -y^2-64z^3
\\
                    16z & -4y& 4x^2-4yz-320xz^2\end{pmatrix}.
\]
The additional equations can be obtained by saturating with respect
to the variable $w$.  The equation of the superisolated singularity
is $\frac14\det M=-\frac12\det N$.
\end{propo}

\section{Degree five}
The first counterexamples to the conjectured formula for $p_g$ of
N\'emethi and Nicolaescu are superisolated singularities with $d=5$
\cite[4.1]{lmn}. The main result of this section is that a splice
quotient with the same resolution graph, if existing, has the
predicted $p_g$.

We first explain how to compute the geometric genus. Consider more
generally a Gorenstein singularity $(X,p)$ such that the exceptional
divisor on the minimal resolution is an irreducible curve $C$ with
arithmetical genus $p_a=g$. The Gorenstein condition implies that
the dual of the normal bundle of $C$ is a line bundle $L$ with
$dL=K_C$ for some integer $d$, where $K_C$ denotes the canonical
sheaf on $C$, so $\deg L = \frac{2g-2}d$. For the canonical cycle on
the minimal resolution we have $K=-(d+1)C$.

We look at the exact sequences
\[
0\lra \sier{}(-(n+1)C) \lra \sier{}(-nC)
\lra \sier{C}(-nC) \lra0
\]
on the minimal resolution. An upper bound for $p_g$ is
$\sum_{n=0}^dh^1(\sier{C}(-nC))= \sum_{n=0}^dh^1(C,nL)$, which by
duality is equal to $\sum_{n=0}^dh^0(C,nL)$. We have equality if the
maps $H^0(\sier{}(-nC))\to H^0(\sier{C}(-nC))$ are surjective. This
is the case if the associated graded ring of the filtration defined
by the vanishing order along $C$ is the ring $\bigoplus H^0(C,nL)$.
This is true by construction for the singularities obtained by
deforming this ring. Splice quotients with the same graph can also
be interpreted in this way.  Okuma has given a formula for $p_g$ for
all splice quotients \cite{okg}.
In  our situation, if $C$ is a rational cuspidal curve of arithmetic
genus 6, and $H$ is the unique line bundle of degree 5 with
$2H=K_C$, then $p_g=7+h^0(C,H)$.

If the irreducible exceptional curve $C$  with $C^2=-d$ has only one
singularity, with only one Puiseux pair $(p,q)$, then a splice type
equation is $x^p-y^q+z^{d+pq}$ with action $\frac 1d(q,p,1)$. The
abstract curve $C$ is given by $x^p=y^q$, where this equation is
considered on a suitable scroll, such that its only singular point
is at the origin of the affine $(x,y)$-chart. If there are two
singular points, with Puiseux pairs $(p,q)$ and $(r,s)$
respectively, then one can again write splice type equations, and
the abstract curve is given by $x^p=y^q$, $z^r=w^s$. If there are
three singular points, the semigroup condition is no longer
satisfied, and one cannot write splice type equations,
\cite[4.3]{lmn}. Indeed, the strict transform of the singular curve
now gives a node in the splice diagram, and the edge weights next to
it are all equal to 1, as in the splice diagram in Proposition
\ref{cuspprop}.

We call an irreducible, locally plane Gorenstein curve hyperelliptic
if the canonical system defines a $2:1$ cover of $\P^1$. The
explicit description above of the curve $C$ on the minimal
resolution  proves the following Proposition.

\begin{propo}
Let the exceptional divisor  of a Gorenstein splice quotient be an
irreducible rational cuspidal curve  $C$ with only singularities of
type $A_k$. Then the curve $C$ is hyperelliptic.
\end{propo}

In an $A_{2k}$-singularity $2k+1$ Weierstra\ss\ points, i.e.,
ramification points of the double cover of $\P^1$,  are absorbed.
The equations for $C$ should be compared with the equations for a
Gorenstein quasi-cone $X(C,L)$, where $C$ is a smooth hyperelliptic
curve and $L$ is a line bundle of degree 1 with $2L$ the $g^1_2$. If
$L=\sier C(P)$ with $P$ a Weierstra\ss\ point, then we get an
equation of type $x^2+y^{2g+1}+z^{4g+2}$. If $L\sim
P_1+\cdots+P_{2k+1}-kg^1_2$, with the $P_i$ Weierstra\ss\ points,
then we get a complete intersection with equations of the type
$x^2=f_{2k+1}(y,w)$, $z^2=g_{2g-2k+1}(y,w)$.

We now specialise to the case $d=5$.

\begin{propo}
Let $(X,p)$ be a splice quotient with exceptional divisor a rational
cuspidal curve  $C$ with $p_a(C)=6$ and $C^2=-5$. Then $p_g(X)$ has
the value predicted by the conjecture of N\'emethi and Nicolaescu.
\end{propo}

Such a singularity exists for every combination of at most two
(higher) cusps. The following table gives the values in question.
\[
\begin{array}{lc@{\qquad\qquad}lc@{\qquad\qquad}lc}
\hline
\mbox{\footnotesize singularities}  & p_g& \mbox{\footnotesize
  singularities}
 & p_g &  \mbox{\footnotesize singularities} & p_g \\[2pt] \hline
\vrule height11.5pt depth3.5pt width0pt
W_{12}  & 10&2E_6    & 7 & A_{10}+A_2 & 8 \\
E_{12}  & 9&E_8+A_4 & 10 & A_8+A_4 & 8\\
A_{12}  & 10&E_6+A_6 & 8 & 2A_6    & 7 \\
\hline
\end{array}
\]
The predicted value for the  cases realisable as plane quintic is
taken from \cite[4.1]{lmn}, the other ones are calculated with the
formulas in \cite[2.3]{lmn}. The geometric genus of the splice
quotients is easily computed. As $p_g=7+h^0(C,5L)$, $\deg L=1$, one
only needs to calculate $h^0(C,5L)$.

Below we comment on several cases. We take the opportunity to give
other singularities with the same resolution graph, not necessarily
superisolated singularities. For some singularities we have
determined the universal abelian cover. To this end we have computed
generators of rings of the type  $\bigoplus H^0(C,nL)$. In most
cases there are many equations between these generators, without
apparent structure.  
It did not
seem worthwhile to compute the corresponding series, or an explicit
example of a universal abelian cover of a superisolated singularity.

\subsection{\boldmath$A_{12}$}
The splice quotient has a hyperelliptic exceptional divisor. Its
equations are described in \cite[4.6]{lmn}. The superisolated
singularity comes from the curve
\[
C\colon\quad z(yz-x^2)^2+2xy^2(yz-x^2)+y^5\;.
\]
Let $H$ be the $g_5^2$ on $C$ and
$L$ the unique line bundle with $5L=H$. We determine generators for
the ring $\bigoplus H^0(C,nL)$ in terms of the parametrisation of
$C$. We start by determining all divisors $D$ of degree 6 such that
$5D$ is cut out by a plane sextic. We find a unique, reduced
divisor, so $\dim H^0(C,6L)=1$, and
$H^0(C,L)=H^0(C,2L)=H^0(C,3L)=0$, as the divisor is not a multiple.
By Riemann-Roch $H^0(C,4L)=0$, and we find $\dim H^0(C,nL)$ for all
other $n$. The following table gives the generators of the ring.
\[
\begin{array}{@{\mbox{degree }}c@{\colon\quad}l}
    5 & t^2s^3+t^5,\; st^4, \; s^5+2s^2t^3, \\[1mm]
    6 & s^6+\frac{12}5s^3t^3+\frac4{75}t^6, \\[1mm]
    7 & s^3t^4+\frac45t^7,
         \; s^7+\frac{14}5s^4t^3-\frac{28}{25}st^6 ,\\[1mm]
    8 & s^4t^4+\frac65st^7,\;  s^6t^2+\frac{11}5s^3t^5+\frac{17}{25}t^8,
        \; s^8+\frac{16}5s^5t^3+\frac{112}{25}s^2t^6, \\[1mm]
    9 & s^3t^6+\frac35t^9, \; s^7t^2+\frac{13}5s^4t^5+\frac{28}{25}st^8,
       \; s^9+\frac{18}5s^6t^3+\frac{88}{25}s^3t^6,
       \; s^5t^4+\frac85s^2t^7\\[1mm]
    11& st^{10},\; s^3t^8+\frac25t^{11},
    \; s^5t^6+\frac75s^2t^9.
 \end{array}
\]
We spare the reader the 104 equations between these generators.

The  general curve of arithmetical genus 6 with an
$A_{12}$-singularity does not have an embedding into the plane. As
example we look at a trigonal curve. Consider the curve
\[
C\colon \quad (w-a^2)^2-2(w-a^2)aw^2+a^4w^3
\]
on $\P^1\times\P^1$ with
parametrisation $w=t^4/(1+2t^3)$, $a=t^2/(1+t^3)$. Its canonical
embedding is
$(1:a:a^2:w:wa:wa^2)$. After homogenising  we have indeed 6 forms in
$s$ and $t$ of degree 10. Let $H$ be the unique line bundle of
degree 5 with $2H=K_C$. It has a unique section, given by
$st^4$. The ring $\bigoplus H^0(C,nH)$ has generators in degree
1, 2 and 3. The corresponding isolated singularity has $p_g=8$.
For the ring $\bigoplus H^0(C,nL)$, with $\deg L=1$, $5L=H$,
one finds that $st^4$ is the generator of lowest degree, so one
needs generators in all degrees from 5 to 13.

\subsection{\boldmath$W_{12}$}
In this case the curve singularity has equation $x^5-y^4-ay^2x^3$.
The results here are similar to the case of an $E_6$-singularity on
a quartic curve.

The quasi-homogeneous case occurs as splice quotient. The equation
is $\xi^5-\eta^4-\zeta^{26}$, and invariants for the $\Z_5$-action
are $x=\xi\zeta$, $y=\eta$ and $z=\zeta^5$, which gives the
superisolated singularity $x^5-y^4z-z^6$.

The non quasi-homogeneous form occurs for the superisolated
singularity $x^5-y^4z-y^2x^3-z^6$. We parameterise the exceptional
curve, and compute the divisors $D$ of degree 6, such that $5D$ is
cut out by a sextic. This time there is a pencil, so the ring
$H^0(C,nL)$, where $\deg L=1$, has one generator in degree 4.

Both singularities occur also on projectively non equivalent curves,
$x^5-y^4z-y^3x^2$ and  $x^5-y^4z-y^2x^3-y^3x^2$. In both cases the
ring $\bigoplus H^0(C,nL)$ has no generators of degree less than 5,
and the number of generators is equal to the number found for the
plane curve with an $A_{12}$.

\subsection{\boldmath$E_{12}$}
There is no plane quintic with a singularity of this type, but we
can form a splice quotient for the graph. It is the quotient of
$\xi^3+\eta^7+\zeta^{26}$ under the action $\frac15(2,3,1)$.

This gives us a clue for writing down rational curves with an
$E_{12}$-singularity. The canonical embedding of the trigonal
rational curve $\xi^3+\eta^7+a\xi\eta^5$  is given by forms of
weighted degree at most 10, where we give the variables $(\xi,\eta)$
the weights $(7,3)$. We set $y_i=\eta^i$, $x_i=\xi\eta^i$. The
equations are in rolling factors format:
\begin{gather*}
\Rank\begin{pmatrix} y_0&y_1 &y_2&x_0\\  y_1&y_2& y_3&
x_1\end{pmatrix} \leq 1 \;,\\ x_0^3+y_3^2y_1+ax_0y_2y_3\;, \\
x_0^2x_1+y_3^2y_2+ax_1y_2y_3\;,\\ x_0x_1^2+y_3^3+ax_1y_3^2\;.
\end{gather*}
For $a=0$ we have the exceptional curve of the splice type
singularity, but for $a\neq0$ the trigonal curve is not invariant
under the $\Z_5$-action. In this case we find that the lowest degree
generator of the ring $\bigoplus H^0(C,nL)$, $\deg L=1$, has degree
6. This means that the ring of $H=5L$ has no generator in degree 1,
but 6 generators in degree 2 and 10 in degree 3. The corresponding
isolated singularity has $p_g=7$.

\subsection{\boldmath$E_8+A_{4}$}
The splice quotient is a superisolated singularity with homogeneous
part $x^3z^2-y^5$. One has also the curve $x(xz-y^2)^2-y^5$, for
which the ring $\bigoplus H^0(C,nL)$ has no generators of degree
less than 5.

\subsection{\boldmath$E_6+A_6$}
The   splice quotient  has high embedding dimension. The
corresponding singular curve is tetragonal. The plane curve
$x(xz+y^2)^2+2y^3(xz+y^2)-y^4z$
leads to $\bigoplus H^0(C,nL)$ without  generators of degree less
than 5.

\subsection{\boldmath$2E_{6}$}
The  exceptional curve of the splice quotient  is a trigonal curve
with bihomogeneous equation $x_1^3y_2^4-x_2^3y_1^4$ on $\P^1\times
\P^1$. As the ring $\bigoplus H^0(C,nL)$ has generators in degree
$3$ and $4$, there are no sections of $5L$, so the  splice quotient
with $E^2=-5$ on the minimal resolution has maximal ideal cycle $2E$
on the minimal resolution.

The most general canonical curve of genus 6 is the complete
intersection of a quintic Del Pezzo surface with a quadric. A plane
model for such a curve is a sextic with 4 double points. An example
with $2$ $E_6$-singularities is
\[
(xz+yz-4xy)^2z^2-4x^3y^3\;,
\]
with a double point in $(1:1:1)$ and three infinitely near double
points at $(0:0:1)$. This curve is still rather special, which
manifests itself in the fact that the ring $\bigoplus H^0(C,nL)$ has
one generator in degree 4.

\subsection{Rational plane quintics with three or four cusps}
The
superisolated singularities with these exceptional curves have all
$p_g=10$, whereas the N\'emethi-Nicolaescu formula gives a lower
value, see the table in \cite{lmn}. In the case $E_6+A_4+A_2$ the
predicted value $8$ is realised by a general curve of arithmetical
genus 6 with these singularities. In the two remaining cases the
predicted value is less than 6. This can never be realised by a
normal surface singularity, as one plus the arithmetical genus of the
exceptional divisor gives a lower bound for the geometric genus of
the singularity.

\section{Generalised superisolated singularities}

The decisive property of the superisolated singularities studied
above is that the exceptional locus of the minimal resolution
consists of one irreducible curve. With the cone over this curve
comes a whole Yomdin series of singularities, whose lowest element
is a superisolated singularity.

A natural generalisation is to look at the series of weighted
homogeneous curves. Let $L$ be an ample line bundle on a singular
curve. The quasi-cone $X(C,L)$ is the singularity with local ring
$\bigoplus H^0(C,nL)$. If this non-isolated singularity is a
hypersurface singularity, we obtain a Yomdin series by adding high
powers of a linear form. If the singular locus has itself singular
branches, this series can be refined. In general, in the non
hypersurface case, there will also be Yomdin type series, as in the
examples for degree four, but there is no easy general formula. To
resolve such singularities one starts with a weighted blow-up. The
singular curve is then the exceptional set.  If this blow-up
resolves the singularity, we speak of a generalised superisolated
singularity.
To obtain a generalised superisolated singularity with rational
homology sphere link from an irreducible singular curve the curve
has to be locally irreducible.

\subsection{Rational curves with an $E_8$-singularity}

A rational curve with an $E_8$-singularity has arithmetical genus
$4$. In its canonical embedding it is a complete intersection of a
quadric and a cubic.
A corresponding superisolated complete intersection  splice quotient
is the quotient of $x^3+y^5+z^{21}$. With variables $d=xz$,
$c=z^3y$, $b=z^6$ and $a=y^2$ we get the equations
$c^2-ab=d^3+a^2c+b^4$. The line bundle $L=\frac12K$ gives rise to a
weighted superisolated singularity $d^3+\alpha^5\beta+\beta^7$,
which is the quotient  of $x^3+y^5+z^{18}$ by a group of order 3.

For both singularities we have a deformation of positive weight,
such that the universal abelian cover is not a deformation of the
splice type equation. We consider $c^2-ab=d^3+a^2c+a^2d+b^4$ and
$d^3+\alpha^5\beta+\alpha^4d+\beta^7$. To study coverings we
parameterise the rational curve $c^2-ab=d^3+a^2c+a^2d=0$ by
\[
(a,b,c,d)=(s^6,(s^2+t^2)^2t^2,-(s^2+t^2)ts^3,ts^5)\;.
\]
Let
$L=\frac16K$. Generators of the ring $\oplus H^0(C,nL)$ are
\[
\begin{array}{@{\mbox{degree }}c@{\colon\qquad}l}
    3 & s^3,\quad s^2t+t^3, \\
    4 & 2s^4+12s^2t^2+9t^4, \\
    5 & 5s^4t+15s^2t^3+9t^5,\quad s^5+3s^3t^2 ,\\
    6 & s^5t, \\
    7 & s^6t, \quad s^7+3s^5t^2.
 \end{array}
\]
There are 20 equations, which we suppress. In particular we obtain
that the singularity $d^3+\alpha^5\beta+\beta^7+t\alpha^4d$ is a
splice quotient only for $t=0$.

The most general rational curve with an $E_8$-singularity lies on a
smooth quadric. Equations are $ad-bc=(a+d)^3+c^2a+\lambda c^2d=0$
with modulus  $\lambda$. Let again $L=\frac16K$. A computation with
the special value $\lambda=0$ shows that $H^0(C,3L)=0$. The ring
$\oplus H^0(C,nL)$ has 23 generators, of degrees $4$ to $9$. The
ring of the double cover corresponding to $3L=\frac12K$ has 4
generators in degree 2 and 6 generators in degree 3. The ideal is
generated by 35 equations.

\subsection{Canonical curves with $g$ cusps}
A parametrisation of a curve with only ordinary cusps is easily
given \cite{st}. Let $\vp_1$, \dots, $\vp_g$ be distinct linear
forms in $s$ and $t$. The map $(\vp_1^2:\dots:\vp_g^2)$ embeds $\P^1$ as
conic in $\P^{g-1}$, which is tangent to the coordinate hyperplanes.
The reciprocal transformation $(z_1:\dots:z_g)\mapsto
(1/z_1:\dots:1/z_g)$ sends it to a $g$-cuspidal curve of degree
$2g-2$.

We consider in particular the case $g=4$. The existence of sections
of roots of the line bundle $K_C$ depends on the modulus of the
curve. First let the  cross ratio of the four cusps on $\P^1$ be
harmonic. We parameterise the canonical curve as follows:
\begin{align*}
  a \;&=\; (t+s)^2(t-is)^2(t+is)^2, \\
  b \;&=\; -(t-s)^2(t-is)^2(t+is)^2, \\
  c \;&=\; i(t-s)^2(t+s)^2(t+is)^2, \\
  d \;&=\; -i(t-s)^2(t+s)^2(t-is)^2,
\end{align*}
and get equations $4ab-(a+b)(c+d)+4cd=abc+abd+acd+bcd=0$. The ring
$\oplus H^0(C,nL)$, where $L=\frac16K$, has $8$ generators, and we
get again 20 equations. The even subring $\oplus H^0(C,2nL)$ gives a
complete intersection. With $x=st$, $y=s^4+t^4$, $z_1=(s^2+t^2)^3$
and $z_2=(s^2-t^2)^3$ we get the equations
$z_1^2-(y+2x^2)^3=z_2^2-(y-2x^2)^3=0$. A weighted superisolated
singularity with exceptional curve a 4-cuspidal rational curve of
self-intersection $-2$ is then given by the equations
$z_1^2-(y+2x^2)^3+x^7=z_2^2-(y-2x^2)^3+x^7=0$. It has Milnor number
65 and $p_g=8$.

For a general cross ratio there are no generators of low degree, and
we get again 23 generators for the universal abelian cover. The ring
$\oplus H^0(C,2nL)$ is also not a complete intersection.

\subsection{A further generalisation}\label{verder}

The methods of this paper can be applied in even more general
situations. We can study rings of the type $\bigoplus H^0(C,nL)$
also for $\Q$-divisors on singular curves. As example we look at the
following graph of \cite[Examples 2]{nw}, which is an example where
the semigroup condition is not satisfied for the splice diagram.
\[
\unitlength=30pt
\def\ci{\circle*{0.26}}
\def\vi{\ci}
\def\mb#1#2{\makebox(0,0)[#1]{$\scriptstyle #2$}}
\begin{picture}(3.5,2.5)(-1.26,-1.25)
\put(-1,-1){\line(1,1){.9}}  \put(1.5,0){\line(-1,0){1.37}}
\put(-1,-1){\vi} \put(-1,-1.2){\mb t{-3}} \put(1.5,0){\ci}
\put(1.5,-.2){\mb t{-7}}
\put(-1,1){\line(1,-1){.9}} \put(-1,1){\ci}
\put(2.5,1.2){\mb b {-3}}
\put(0,0){\circle{0.26}} \put(0,-.2){\mb t{-1}}
\put(2.5,1){\line(-1,-1){1}}   \put(2.5,1){\vi}
\put(2.5,-1.2){\mb t{-3}} \put(2.5,-1){\line(-1,1){1}}
\put(2.5,-1){\vi}
\end{picture}
\qquad
\def\ci{\circle{0.2}}
\begin{picture}(5,1)(-4,-1.25)
\put(0,0){\ci}
\put(.1,.04){\line(3,1){1.3}} \put(.3,.22){\mb b3}
\put(.1,-.04){\line(3,-1){1.3}}
\put(.3,-.22){\mb t3} \put(1.5,.505){\ci} \put(1.5,-.505){\ci}
\put(-.1,0){\line(-1,0){1.3}}  \put(-1.2,.1){\mb b{57}}
\put(-.3,.1){\mb b1} \put(-1.5,0){\ci}
\put(-1.6,.04){\line(-3,1){1.3}} \put(-1.8,.22){\mb b2}
\put(-1.6,-.04){\line(-3,-1){1.3}} \put(-1.8,-.22){\mb t3}
\put(-3,.505){\ci} \put(-3,-.505){\ci}
\end{picture}
\]
A Gorenstein singularity with this resolution graph exists. To find
it we first look at the minimal resolution graph (as described in
the first section):
\[
\unitlength=30pt
\def\vi{\circle*{0.26}}
\def\mbt#1{\makebox(0,0)[b]{$\scriptstyle #1$}}
\begin{picture}(2,.5)(-.26,-.25)
\put(0,0){\line(1,0){.87}}  \put(2,0){\line(-1,0){.87}}
\put(1,-.2){\makebox(0,0)[t]{$\scriptstyle [1]$}}
\put(1,0){\circle{0.26}} \put(0,0){\vi} \put(2,-0){\vi}
\put(0,.2){\mbt{-3}} \put(1,.2){\mbt{-1}} \put(2,.2){\mbt{-3}}
\end{picture}
\]
The weight below the central vertex is the arithmetic genus of the
central cuspidal curve $E_0$. This is a star-shaped graph. Consider
the $\Q$-divisor $D=P-\frac13Q_1-\frac13Q_2$ on  $E_0$, where $Q_1$
and $Q_2$ are the intersection points with the other components
$E_1$ and $E_2$, and $P$ is determined by the normal bundle. We have
$K=-5E_0-2E_1-2E_2$, so the Gorenstein condition is $4P-2Q_1-2Q_2=0$
in $\Pic(E_0)$. We satisfy it by taking $2P=Q_1+Q_2$. We compute the
graded ring $\bigoplus H^0(E_0,\lfloor nD\rfloor)$ for $E_0$ a
cuspidal rational curve. We find the hypersurface singularity
$z^2=(x^2-y^3)^3$. An isolated singularity with the original
resolution graph is obtained by adding  generic terms of higher
weight. An example is $z^2=(x^2-y^3)^3+xy^8$. The branch curve is
irreducible and has Puiseux pairs $(3,2)$ and $(10,3)$.

The topology of  the universal abelian cover  is easily computed.
The minimal resolution has as exceptional divisor a 3-cuspidal
rational curve of self-intersection $-1$. Such a singularity was
studied above in section \ref{cusp}. We use the notation and
equations introduced there. A suitable $\Z_3$-action on the
3-cuspidal curve $C$ is a cyclic permutation of the cusps, i.e., a
cyclic permutation of the coordinates on $\P^2$.  The polynomials
$Q$ and $T$ are invariant under all permutations. The polynomial
$27(x^2y+y^2z+z^2x-x^2z-y^2x-z^2y)$ is a skew invariant, which gives
a section $\psi$ of $3K-D\equiv 9(K-D)$. The cyclic permutation
induces an action on the coordinates $(s,u,y,z,v,w)$ of the
equations (\ref{dca},b). As the variables $s$ and $u$ correspond to
$Q$ and $T$, they are invariant. Furthermore $y\mapsto z  \mapsto
s^2-y-z$, and $v$, $w$ and $us$ are involved in more complicated
formulas. The ring of invariants of the action on $\bigoplus
H^0(C,nL)$ is generated by $s$, $u$ and $\psi$, and there is one
equation $$27\psi^2+(4s^3+u^2)^3=0\;.$$ An explicit isolated
singularity with $C$ as exceptional divisor with $C^2=-1$ is given
in Proposition \ref{cuspprop}. Its equations are written explicitly
in the introduction. As they are obtained by the invariant
perturbation $\vp_1=us^2$, the same $\Z_3$-action  on the
coordinates $(s,u,y,z,v,w)$ acts on the singularity, and we find the
quotient by expressing $\psi^2$ in the local ring of the singularity
in terms of $s$ and $u$. We get a hypersurface singularity, which is
equisingular with $z^2=(x^2-y^3)^3+xy^8$, but has a more complicated
equation.
\begin{propo}
\label{verderprop} The isolated singularity of Proposition
\ref{cuspprop} is the universal abelian cover (of order 3) of the
double point
\[
27\psi^2+(4s^3+u^2)^3+2u^5s^2-20u^3s^5-4us^8+u^4s^4\;,
\]
whose branch curve is irreducible, with Puiseux pairs $(3,2)$ and
$(10,3)$.
\end{propo}

\section{Discussion}
Superisolated singularities allow us to import the theory of
projective curves into the subject of surface singularities. Our
singularities with rational homology sphere links do not behave as
rational curves, but have the character of curves of their
arithmetical genus. The conjectures of Neumann and Wahl ask for the
existence of special divisors (those cut out on the complete
intersection by the coordinate hyperplanes). For smooth curves this
is a complicated problem; the only thing we know for sure is that
those of low degree do not exist on curves with generic moduli.

The idea behind the Casson Invariant Conjecture of Neumann and Wahl
is that the Milnor fibre of a complete intersection is a natural
four-manifold attached to its boundary, the link of the singularity,
whose signature in the integral homology sphere case computes the
Casson invariant of the link exactly. Starting from a graph the
construction of Neumann and Wahl is certainly very natural. Our
computations for graphs of superisolated singularities with $d=5$
support the following conjecture.

\begin{conj} Among the $\Q$-Gorenstein singularities with
rational homology link the splice quotients have geometric genus as
predicted by the Seiberg-Witten Invariant Conjecture of N\'emethi
and Nicolaescu.
\end{conj}

This conjecture has now been proved by N\'emethi and Okuma, first in
the  integral homology sphere case \cite{no}, where it reduces to
the Casson Invariant Conjecture, and later in general. More
generally, Braun and N\'emethi prove an equivariant version
\cite{bn}. An important ingredient is Okuma's $p_g$-formula
\cite{okg}.

 In all examples we computed, where the universal abelian cover
is not of splice type, it is not a complete intersection. The same
is true for our Gorenstein singularities with integral homology
sphere link. The following conjecture of Neumann and Wahl is still
open.

\begin{conj} A complete intersection surface singularity with
integral homology sphere link has equations of splice type.
\end{conj}

On the other hand, for all graphs  we studied we found a singularity
which is not a complete intersection. In fact, we conjecture that
this is the general behaviour. If the exceptional divisor is a
reduced and irreducible curve, one takes this curve general in
moduli (whatever this means, as it is no so clear whether a sensible
moduli space exists for curves with given singularities). Of course
this does not apply to rational and elliptic singularities. We also
exclude the case that the curve is hyperelliptic.

\begin{conj} The `general'  $\Q$-Gorenstein singularity
with given rational homology sphere link is not a splice quotient,
if the fundamental cycle has arithmetic genus at least 3. Its
universal abelian cover is not a complete intersection.
\end{conj}

\end{document}